\documentclass{amsart}

\newtheorem{satz}{Theorem}
\newtheorem{lemma}{Lemma}
\newtheorem{kor}{Corollary}
\newtheorem{prop}{Proposition}

\theoremstyle{definition}
\newtheorem{defn}{Definition}

\theoremstyle{remark}
\newtheorem{rem}{Remark}

\numberwithin{equation}{section}

\newcommand{\C}{{\mathbb C}}
\newcommand{\R}{{\mathbb R}}
\newcommand{\Z}{{\mathbb Z}}

\newcommand{\N}{{\mathbb N}}
\newcommand{\CP}{\C P}
\newcommand{\cpq}{\overline{\C P^2}}

\DeclareMathOperator{\U}{U}
\DeclareMathOperator{\SO}{SO}
\DeclareMathOperator{\GL}{GL}
\DeclareMathOperator{\mm}{min}

\DeclareMathOperator{\Tor}{Tor}

\begin{document}

\title{Embedded surfaces and almost complex structures}
\author{Christian Bohr}
\address{Mathematisches Institut der LMU \\
Theresienstr. 39 \\ 80333 M\"unchen \\
Germany}
\email{bohr@rz.mathematik.uni-muenchen.de}
\date{October 20, 1998}

\begin{abstract}
In this paper, we prove necessary and sufficient conditions for a 
smooth surface in a smooth 
4-manifold $X$ 
to be pseudoholomorphic with respect to an almost complex structure on $X$. 
In particular, 
this provides
a systematic approach to the construction of pseudoholomorphic curves that 
do not 
minimize the genus in
their homology class.
\end{abstract}

\maketitle

\section{Introduction and summary of results}

Let $X$ be a closed differentiable and connected 4-manifold 
with an orientation  and $\Sigma \subset X$ a connected
oriented surface. An {\bf almost complex structure $\mathbf
J$} on
$X$ is a reduction of the structure group $\GL^+(4)$ of 
$TX$ to the group
$\GL(2,\C)$, 
or, in other words, a section $J$ of $End(TX)$ with $J^2 = -1$ that
preserves the orientation, so that $TX$ carries the structure
of a complex vector bundle.
The surface $\Sigma$ is called a {\bf pseudoholomorphic curve} if the
tangent bundle of $\Sigma$ is preserved by $J$ (note that in this case,
the almost complex structure on $X$ induces a complex structure on $\Sigma$).
The question that shall be treated on the following pages is : Given a
surface $\Sigma$, is there an  almost complex structure $J$ on $X$
such that $\Sigma$
is a pseudoholomorphic curve with respect to $J$?

\bigskip

First, recall that an almost complex structure $J$ has a first Chern class
$c_1(J) \in H^2(X;\Z)$ (this is just the first Chern class of $TX$ considered
as a complex vector bundle) and that this class has the properties 
\begin{enumerate}
\item $c_1(J)^2 = 2 \chi + 3 \tau$
\item $c_1(J) \equiv w_2 \mod 2$
\end{enumerate}
where $\chi$ denotes the Euler characteristic  and $\tau$ the signature
of the intersection form of $X$. A class with property 2 is called a {\bf
characteristic class} on $X$. If the homology of $X$ does not contain 
2-torsion, then 
these classes can be characterized in terms of
the intersection form $Q$ of $X$: a class $c \in H^2(X;\Z)$ is characteristic
if and only if $Q(x,x) \equiv Q(x,c) \mod 2$ for all $x \in H^2(X;\Z)$.
If there is 2-torsion, one part of this statement is still true: if $c$ is
characteristic, then $Q(x,c) \equiv Q(x,x) \mod 2$ for every $x$. Conversely,
if a class $c$ fulfills $Q(x,c) \equiv Q(x,x) \mod 2$ for all $x$, then there
is a   torsion class $a$ such that $c+a \equiv w_2(X) \mod 2$.
It is a classical result of Whitney that there are characteristic classes on any
4-manifold (\cite{W}).

Furthermore, it is well known that in turn
every class in $H^2(X;\Z)$ fulfilling the conditions 1,2 can be realized as 
the
first Chern class of an almost complex structure. So there is an almost complex
structure on $X$ if an only if there is a class in $H^2(X;\Z)$ that 
fulfills the
conditions above (this is a result of Wu, see \cite{HH}), in fact every
such class is the first Chern class of an almost complex structure.
Consideration of the
intersection form easily leads to the conclusion 
that if the intersection form of $X$ is indefinite, there is an almost complex
structure on $X$ if and only if $b_1+b_2^+$ is odd, where $b_2^+$ denotes the
maximal dimension of a subspace of $H^2(X,\R)$ on which the intersection form 
is positive definite.

\begin{defn} Let $G$ be a finitely generated abelian group and $g \in G$.
Let $T \subset G$ be the torsion subgroup of $G$.
\begin{enumerate}
\item If $G$ is free abelian and $g \neq 0$, the {\bf divisibility} of $g$ 
in $G$ is 
defined to be the largest positive integer
$d$ with the property that there is an $x \in G$ with $g=dx$. 
The divisibility of $0 \in G$ is
defined to be zero.
\item For arbitrary $G$, the divisibility $d(g)$ of $g \in G$ is the 
divisibility of $g$
(more precisely, the residue class of $g$) in the free abelian group $G/T$. 
The divisibility is defined to be zero if and only if $g \in T$.
\end{enumerate}
\end{defn}

\begin{rem}\label{rem1.1}
 Clearly the image of the homomorphism $Hom(G;\Z) \rightarrow \Z$,
given by evaluation on $g$ , is just $d(g) \Z$. From this, we see that the 
divisibility of $k \cdot g$ for $g \in G$, $k \in \Z$ is $\pm k$ times the
divisibility of $g$.
\end{rem}

\begin{defn} Let $(\Gamma, Q)$ be a lattice (i.e. $\Gamma$ is a free abelian
group of finite rank and $Q$ a unimodular symmetric bilinear form on $\Gamma$).
For $\gamma \in \Gamma$ with divisibility $d=d(\gamma)$ 
define  $k(\gamma) \in \Z_d$ as follows: Choose a
characteristic class $c \in \Gamma$, i.e. $Q(c,x) \equiv Q(x,x) \mod 2$ for
every $x \in \Gamma$,  and set
$$
k(\gamma) := 1 + \frac{1}{2}(Q(\gamma,\gamma) - Q(c,\gamma)) \mod d
$$
\end{defn}

The residue class $k(\gamma)$ is independent of the choice of $c$: 
if $c'$ is another
characteristic class, then $c$ and $c'$ differ by a multiple of $2$, so
the terms in the bracket differ by a multiple of $2d$, 
according to remark \ref{rem1.1}, and this does not affect $k(\gamma)$. 
If $\gamma=0$, then
$k(\gamma)=1 \in \Z$.

\begin{defn} For a closed connected and oriented surface $\Sigma \subset X$
let $k(\Sigma)=k([\Sigma])$ with respect to the lattice defined by the
intersection form on the free group $H_2(X;\Z) / \Tor H_2(X;\Z)$, 
where $[\Sigma]$ denotes the
homology class of $\Sigma$.
\end{defn}

Since cup products are not altered by adding torsion classes
to one of the factors, we could as well have defined $k([\Sigma])$ by
$$
k([\Sigma]) = 1 + \frac{1}{2}(\Sigma \cdot \Sigma - c \cdot \Sigma) \mod d
$$ 
where $d$ denotes the divisibility of $[\Sigma]$ in $H_2(X;\Z)$ and $c$
is any characteristic class on $X$, i.e. $c \equiv w_2(X) \mod 2$. We will
use the notation $d(\Sigma)$ for the divisibility of the class $[\Sigma]$.

If there is an almost complex structure $J$ turning $\Sigma$ into a
pseudoholomorphic curve, then the adjunction formula
$$
g(\Sigma) = 1 + \frac{1}{2} (\Sigma \cdot \Sigma - c_1(J) \cdot \Sigma)
$$
holds, and $c_1(J)$ is characteristic, so  we have the congruence
$$
g(\Sigma)
\equiv k(\Sigma) \mod d.
$$
It turns out that this necessary condition is in
fact sufficient for the existence of such a $J$ if the intersection form of 
$X$ is strictly 
indefinite (i.e. 
$\mm \{b_2^+,b_2^-\} \geq 2$):

\begin{satz}\label{thm1}
Let $X$ be a connected oriented closed and differentiable 4-manifold
and $\Sigma \subset X$ a closed connected
and oriented surface with divisibility $d$. Suppose $\mm \{ b_2^+,b_2^- \}
\geq 2$ and $b_1+b_2^+ \equiv 1 \mod 2$. Then
there is an almost complex structure $J$ on $X$ such that the
surface $\Sigma$ is pseudoholomorphic with respect to $J$ if and only if
$g(\Sigma) \equiv k(\Sigma) \mod d$.
\end{satz}

Note that this condition is in particular fulfilled when the class of
$\Sigma$ has divisibility one, 
so any such surface is pseudoholomorphic with
respect to an almost complex structure on $X$.
In addition, if $[\Sigma]$ is not a torsion class, 
this condition is ``cyclic'' :  
we can attach handles that do not change the homology class of $\Sigma$ - and
hence do not alter $k(\Sigma)$ - but
raise the genus until the condition of the Theorem is fulfilled. In this way
we even can produce a surface $\Sigma'$
homologous to $\Sigma$ that is pseudoholomorphic with respect to an almost
complex structure on $X$, but whose genus is arbitrarily large:

\begin{kor} Let $X$ be a closed connected and oriented smooth 4-manifold 
 with $b_1+b_2^+ \equiv 1 \mod 2$ and $\mm \{
b_2^+,b_2^- \} \geq 2$, and let $\Sigma$ be a surface in $X$ such 
that $[\Sigma]$ is not a torsion class. Let 
$m \in \N$ be any natural number. 
Then there is an almost complex structure $J$ on $X$ and a pseudoholomorphic
curve $\Sigma'$ homologous to $\Sigma$ with $g(\Sigma') \geq g(\Sigma)
+ m$.
\end{kor}

This Corollary 
 provides a large number of pseudoholomorphic curves
that do not minimize the genus in their homology class. Other examples for 
this have 
been given by Kotschick 
(unpublished) and in a paper by 
Mikhalkin
(\cite{Mi}). 
Although the case $X=\C P^ 2$ is not covered by the Corollary, this  can 
also occur on $\CP^2$, 
an example for this is  
the curve in the statement of Proposition \ref{prop3}.

The next three Propositions show that the condition 
$\min\{b_2^+,b_2^-\} \geq 2$
is really necessary, if it is dropped, the Theorem is no longer true:

\begin{prop}\label{prop1} 
If $X$ is a rational complex surface, there is a surface $\Sigma$ in $X$
with $g(\Sigma) \equiv k(\Sigma)\mod d(\Sigma)$ 
that is not pseudoholomorphic with
respect to any almost complex structure on $X$.
\end{prop}

\begin{prop}\label{prop2} Let $X$ be a 4-manifold with definite 
intersection form. 
Then there is a surface $\Sigma$
in $X$ with $k(\Sigma) \equiv g(\Sigma) \mod d(\Sigma)$ 
that is not pseudoholomorphic
with respect to any almost complex structure on $X$.
\end{prop}

\begin{prop}\label{propcountex} 
Let $Q$ be a unimodular symmetric bilinear form over $\Z$ 
fulfilling $\min\{b^+,b^-\} \leq 1$
that can be realized as the
intersection form of a smooth 4-manifold. 
In the case that $Q$ is indefinite and even
assume that the signature of $Q$ is non-negative.
Then there is a closed oriented 4-manifold $X$ having $Q$ as 
intersection form and a closed oriented
and connected surface $\Sigma \subset X$ such that
$b_1(X) + b_2^+(X) \equiv 1 \mod 2$ and $g(\Sigma) \equiv 
k(\Sigma) \mod d(\Sigma)$, but $\Sigma$ is not pseudoholomorphic 
with respect to
any almost complex structure on $X$.
\end{prop}

Note that all odd forms and all even intersection forms of smooth 4-manifolds 
that have  no 2-torsion in 
their homology (\cite{D}) or are spin (\cite{Fu}) are covered by this
Proposition.

Finally, there is a simple example for the case $X = \CP^2$. 
The class $-1 \in H_2(\CP^2;\Z)$ can be
represented by a sphere - just take the complex line with the orientation
reversed -, hence the minimal genus for this class is $0$. The following
Proposition therefore  provides
another example that a pseudoholomorphic curve does not always minimize
the genus in its homology class:

\begin{prop}\label{prop3} 
There is a surface with genus 3, representing minus the generator 
of $H_2(\CP^2;\Z)$, that is pseudoholomorphic with respect to an 
almost complex
structure homotopic to the canonical one.
\end{prop}

\section{Proofs of Theorem \ref{thm1} and Proposition \ref{prop3}}

For the proofs, we need two Lemmas, the first of them being a topological
Lemma, whereas the second one is purely algebraic:

\begin{lemma}\label{lemmahom} If there is a class $c \in H^2(X;\Z)$ with 
the following
properties
\begin{enumerate}
\item $c^2=2\chi+3\tau$
\item $c \equiv w_2 \mod 2$
\item $<c,[\Sigma]> = 2 -2g(\Sigma) + \Sigma \cdot \Sigma$,
\end{enumerate}
then there is an almost complex structure $J$ such that 
the surface $\Sigma$ is 
pseudoholomorphic with respect to $J$ and $c_1(J)=c$ (here $\tau$ denotes
the signature of $X$).
\end{lemma}

\begin{proof}
By the result of Wu mentioned earlier (see \cite{HH}), there is an 
almost complex structure 
$J_0$ on
$X$ with 
$c_1(J_0)=c$. 
Introduce a Riemannian metric $g$ on $X$ compatible with $J_0$, i.e. the 
endomorphism
$J$ is isometric on the fibers of $TX$ with respect to $g$. Then the
almost complex structures compatible with $g$ can be identified with
the reductions of the structure group $\SO(4)$ to $\U(2)$, i.e. with sections
in the bundle $\Theta:=P_{\SO(4)} / \U(2)$ having fiber $\SO(4) / \U(2) = 
S^ 2$.
For the restriction of the tangent bundle to $\Sigma$, we have a 
decomposition $TX | \Sigma = N \oplus T\Sigma$, where $N$ denotes the
normal bundle of $\Sigma$. By introducing metrics on these two bundles,
their structure group can be reduced to $\SO(2)$. Since $\SO(2) \times \SO(2)
 \subset \U(2)$, we have an almost complex structure on $TX | \Sigma$ turning
this decomposition into a direct sum of complex vector bundles. This almost
complex structure 
can be extended to an almost complex structure $J$ on 
the disk bundle $DN$ (that is identified with a tubular
neighborhood of $\Sigma$). 
Clearly, $\Sigma$ is a pseudoholomorphic curve in $DN$ with
respect to $J$. We now have to show that $J$ can be extended over $X$
to an almost complex structure homotopic to $J_0$ as a section of $\Theta$,
then $c_1(J)=c_1(J_0)=c$, and the Lemma is proved.

The second cohomology $H^2(DN;\Z)$ is $\Z$, generated by the fundamental
class $[\Sigma]$ (more exactly, by its pullback via the projection $DN 
\rightarrow \Sigma$). Let $s$ respectively $s_0$ denote the sections
of $\Theta$ on $DN$ given by $J$ and $J_0$.
Note that $J_0$ defines an extension of $s_0$ to $X$.
Let $c_1 \in H^ 2(DN;\Z)$ denote the first
Chern class of $J$. By definition of $J$, we have a decomposition 
$(TX | \Sigma, J) = N \oplus T\Sigma$ of complex vector bundles. Taking
the first Chern class on both sides yields the adjunction equality 
$<c_1, \Sigma>= 
2-2g + \Sigma \cdot \Sigma$. But by assumption 3, the same is true 
for $c=c_1(J_0)$, 
hence $c_1=c$ in $H^ 2(DN;\Z)$. A short calculation, using the
exact homotopy sequence of the fibration
$$
S^ 2 \rightarrow B\U(2) \rightarrow B\SO(4),
$$
yields that $\pi_2(S^ 2) \rightarrow \pi_2(B\U(2))$ is the multiplication by
2, and this shows that for the primary difference $p \in H^2(DN;\Z)$ 
between $s$ and $s_0$ 
as sections $DN \rightarrow \Theta$, we have the equality
$2p = c_1 - c=0$. Since the homology of $DN$ is torsion free, this implies
$p=0$, and since $H^ 3(DN;\Z)=H^ 4(DN;\Z)=0$, there are no higher obstructions,
hence the sections $s$ and $s_0$ are homotopic on $DN$. Using the 
homotopy extension property we can conclude that there is an extension
of $s$ to $X$ homotopic to $s_0$, and this proves the assertion. 
\end{proof}

\begin{lemma}\label{lemmaalg} Let $(\Gamma,Q)$ be a lattice with 
$\min\{b^+,b^-\} \geq 2$, let $\gamma \in \Gamma$ be a vector with
divisibility $d$, $h$ an integer with $h \equiv \tau(Q) \mod 8$, where
$\tau(Q)$ denotes the signature of $Q$, and $g$ be a natural number with
$g \equiv k(\gamma) \mod d$.
Then there is a $c \in \Gamma$ with
\begin{enumerate}
\item $c$ is characteristic, i.e. $Q(c,x) \equiv Q(x,x) \mod 2$ for every
$x \in \Gamma$,
\item $Q(c,c) = h$ and
\item $Q(c,\gamma) = 2 -2g + Q(\gamma,\gamma)$.
\end{enumerate}
\end{lemma}

\begin{proof} According to the classification Theorem 
of Hasse-Minkowski 
(see \cite{MH}), we can choose a basis $(e_1, \dots , e_n)$ such
that with respect to this basis, $Q$ is described by the matrix
$$
\begin{pmatrix} \framebox {$\begin{matrix} & & \\ & H & \\ & & 
\end{matrix}$} & & \\ 
& \framebox {$\begin{matrix} & & \\ & H & \\ & & \end{matrix}$} & \\
& & \framebox {$\begin{matrix} & & \\ & A & \\ & & \end{matrix} $}
\end{pmatrix}
$$
where $H = \begin{pmatrix} 0 & 1 \\  1 & 0 \end{pmatrix}$ denotes the 
hyperbolic form, and $A$ is diagonal if $Q$ is odd, or of the type
$mE_8$ with some $m \in \Z$ in the case that $Q$ is even.
If $\gamma=0$, the condition on $g$ reads $g=1$, and any characteristic 
$c$ with $Q(c,c)=h$ will do the job (it is
easy to see that such a $c$ exists). Now assume $\gamma \neq 0$. 
Let $\gamma=dp$ with a $p \in \Gamma$ having divisibility
one. 

\noindent {\bf Case 1: } $Q$ is even. Then $p$ must be ordinary (i.e. not
characteristic), because $Q$ is unimodular and - according to the 
characterisation of 
the divisibility given in Remark \ref{rem1.1} - therefore
there is an $x \in \Gamma$ with $Q(x,p)=1$. Using a result of Wall
(\cite{W1}) concerning the group of automorphisms of $Q$, we can assume 
$p=(k,1,0, \dots, 0)$ with some $k \in \Z$ (Wall's Theorem asserts that
there is an automorphism that maps $p$ to some vector of this type, 
after a change of the basis, we can assume that $p$ has this special form).
Let $c_0 \in <e_3, \dots, e_n>$ be some characteristic vector with
$Q(c_0,c_0)=h$ (it is easy to see that such a $c_0$ exists, 
using the Hasse-Minkoswki classification applied to $H \oplus A$).The 
assumption on $g$ implies
that the difference between $Q(c_0,\gamma)$ and $2-2g+Q(\gamma,\gamma)$ is
a multiple of $2d$, say $2da$ with $a \in \Z$. Let $c=c_0 + 2ae_1$.
Then $Q(c,c)=Q(c_0,c_0)=h$ and $Q(c,\gamma)=Q(c_0,\gamma)+2ad=
2_2g+Q(\gamma,\gamma)$.

\noindent {\bf Case 2:} $Q$ is odd:

\noindent a) $p$ is ordinary: Then, again using the result of Wall, we
can assume that $p=(k,1, p')$ with $p' \in <e_3, \dots, e_n>$, and the
same arguments as in Case 1 apply.

\noindent b) $p$ is characteristic: Since $Q$ is odd, the same must be 
true for $A$, 
in particular, $n \geq 5$. In this case, the standard form 
for
$p$ is $p=(0,0,2k,2,1, \dots, 1)$ with some $k \in \Z$, because this vector
is characteristic, has divisibility one and square $8k+\tau(Q)$, 
and $Q(p,p) \equiv \tau(Q) \mod 8$. 
Choose any characteristic vector $c_0 \in <e_3, \dots , e_n>$. 
Then the difference $Q(c_0,c_0) - (2-2g+Q(\gamma,\gamma))$ is divisible by
$2d$ (this follows from the assumption on g and the definition of 
$k(\gamma)$). 
Therefore we can choose $a \in \Z$ such that $c_1:=c_0+2ae_5$ has the
properties 1 and 3 (observe $Q(e_5,p)=\pm1$). The difference between
$Q(c_1,c_1)$ and $h$ is now a multiple of 8, say $8b$, $b \in \Z$, and
therefore $c=c_1+2be_1+2e_2$ fulfills all the three conditions.
\end{proof}

\begin{proof} of Proposition \ref{prop3}:
Let $J$ denote the standard almost complex structure on $\CP^2$.
with Chern class $c_1(J)=3$. 
Let $\Sigma'$ denote the complex line $\CP^1$ in $\CP^2$ with the
orientation reversed, hence $[\Sigma']=-1$. By attaching
handles, we can construct a surface $\Sigma \subset
 \CP^2$ with genus 3 homologous
to $\Sigma'$. For this surface, we have
$$
<c_1(J),[\Sigma]> = -3 = 2-2g(\Sigma) + \Sigma \cdot \Sigma,
$$
and the assertion follows using the homotopy argument as in the
proof of Lemma  \ref{lemmahom}.
\end{proof}

\begin{proof} of Theorem \ref{thm1}: First assume that $\Sigma$
is pseudoholomorphic with respect to $J$. Then we have the adjunction
equality $c_1(J) \cdot \Sigma = 2 -2 g + \Sigma \cdot \Sigma$, and this
implies $g(\Sigma) \equiv k(\Sigma) \mod d$. For the converse, let $\Sigma$
fulfill $g(\Sigma)=k(\Sigma) \mod d$. 
Let $\Gamma=H^ 2(X;\Z) / \Tor H^ 2(X;\Z)$ and $Q$ denote the form on $\Gamma$
defined by the intersection form. Let $\gamma \in \Gamma$ be the residue
class of $[\Sigma]$ and $h:=2 \chi + 3 \tau$.
A short calculation shows that the condition $b_1+b_2^ + \equiv 1 \mod 2$ 
implies that $\chi + \tau $ is divisible by $4$ and $h \equiv \tau(Q) \mod 8$.
According to Lemma \ref{lemmaalg}, there is a $c' \in \Gamma$ with
$Q(c',c') = h, Q(c',\gamma) = 2-2g + Q(\gamma,\gamma)$ and $Q(c',x) 
\equiv Q(x,x) \mod 2$ for all $x \in \Gamma$. Choose a lift $c \in 
H^ 2(X;\Z)$ of
$c'$ such that $c \equiv w_2(X) \mod 2$. Then $c$ fulfills the conditions
of Lemma \ref{lemmahom}, and the assertion of the Theorem follows.
\end{proof}

\section {\bf Proof of Propositions \ref{prop1},\ref{prop2} and 
\ref{propcountex}}

\begin{proof} of Proposition \ref{prop1}: 
A rational surface is diffeomorphic
to $S^2 \times S^2$ or to $\CP^2 \# k\cpq$. First, consider the case $X = S^2
\times S^2$. Let $\Delta$ denote the diagonal sphere in $S^2 \times S^2$. 
Whenever $c=(x,y) \in H^2(X;\Z)=\Z \times \Z$ is the Chern class of an almost
complex structure, we have $c^2=2xy=8$, and $x,y$ are even, therefore 
we have $c=(2,2)$ or $c=(-2,-2)$. 
Now choose a surface $\Sigma$ homologous
to $\Delta$. Observe that, since $\Delta$ is a sphere, one can
construct such surfaces $\Sigma$ of any genus by attaching nullhomologous
handles to $\Delta$. 
Then, for any almost complex structure $J$ on $X$, we have
$c_1(J) \cdot \Sigma=\pm 4$ and $\Sigma \cdot \Sigma=2$. 
If we choose $\Sigma$ to have genus
0 or 4, we see that there is an almost complex structure on $X$ that turns
$\Sigma$ into a pseudoholomorphic curve, but if we choose a surface $\Sigma
\sim \Delta$ with genus 1, then there is no almost complex structure on $X$
such that $\Sigma$ is pseudoholomorphic. But on the other hand, the
divisibility of $\Sigma$ clearly is one, so the equality $g(\Sigma) \equiv
k(\Sigma) \mod d$ is fulfilled for every value of $g(\Sigma)$.

Now consider the case $X=\CP^2$. 
Clearly, only the classes $3$ and $-3$ occur as Chern classes of almost complex
structures on $X$. 
If we construct a surface $\Sigma$ of genus one, representing the generator
of $H^2(X;\Z)$, by attaching a handle to $\CP^1$, then, for any almost complex
structure $J$, $c_1(J) \cdot \Sigma = \pm 3$, $\Sigma \cdot \Sigma = 1$, so
the adjunction equality will not hold for $J$. Again, we have $d=1$,
and this provides the required example.

Now we turn to the case $X=\CP^2 \# \cpq$. Let $\Sigma$ be a surface of
genus one, representing the class $(1,0) \in H^2(X;\Z)=\Z \times \Z$. 
Then the divisibility is 1, and we have to prove that there is
no almost complex structure such that $\Sigma$ is pseudoholomorphic.
If $J$ is an almost complex structure on $X$ with first Chern class
$c=(x,y)$, then $x^2 - y^2 = 8$. If $\Sigma$ would be pseudoholomorphic
with respect to $J$, this would imply $c \cdot \Sigma = x = 1$, hence
$1 - y^2 = 8$, contradiction. This proves that  there is no such $J$.

Finally, to settle the case $X=\CP^2 \# k\cpq$ and $k \geq 2$, note that
$X$ is diffeomorphic to $(S^2 \times S^2)\# (k-1) \cpq$ (\cite{W2}). 
 Let again $\Delta
\subset X$ denote the sphere coming from the diagonal embedding in 
$S^2 \times
S^2$, choose a surface $\Sigma$ representing the
same homology class and let $g$ denote its genus. With respect to the
basis of $H^2(X;\Z)$ coming from a diffeomorphism $X \cong (S^2 \times
S^2)\# (k-1) \cpq$, we have $[\Sigma]=(1,1,0, \dots , 0)$, and if $J$ is
an almost complex structure with Chern class $c=(x,y,a)$, with 
$a \in H^2((k-1)\cpq)$ and $x,y \in \Z$, 
then $c^2 = 2xy + a \cdot a = 9 - k$ (note that $a \cdot a \leq 0$, here
the dot denotes the cup product in the cohomology of $(k-1)\cpq$).
Now $a \cdot a  \leq -k+1$, since $a \equiv w_2((k-1)\cpq) \mod 2$ and the 
intersection form of $(k-1)\cpq$ is standard,
and we can conclude  
$$
2xy = 9 - k - a \cdot a \geq 8.
$$
Therefore, as in the example $X=S^2 \times S^2$, we have $xy \geq 4$, and
$x,y \equiv 0 \mod 2$.
Now suppose that $\Sigma$ is pseudoholomorphic with respect to $J$. Then
we have the adjunction equality $c \cdot \Sigma = x + y = 4 -2 g$. 
Together with $xy \geq 4$, this implies $g=0$ or $g \geq 4$. But we can 
construct $\Sigma$ by attaching one handle at $\Delta$ and therefore realize
$g=1$. Hence this surface is not pseudoholomorphic with respect to any
almost complex structure on $X$, and this completes the proof. 
\end{proof}

Note that the last part of the proof can be applied to every $X$ of
the type $(S^2 \times S^2) \# N$, where $N$ has no 2-torsion
in its homology and negative definite intersection
form (which must be standard, according to Donaldson).

\begin{proof} of Proposition \ref{prop2}: 
The intersection form $Q$  can be considered as a non-degenerate 
symmetric bilinear form on
the real cohomology $H^2(X;\R)$, where the free part of the integral 
cohomology is
lying as a lattice in this real vector space. Choose a class $\gamma 
\in H^2(X;\Z)$
with self-intersection $s=\gamma \cdot \gamma \neq 0$ and divisibility one. 
First suppose that $Q$ is positive
definite. Then, for every class $c \in H^2(X;\Z)$, we have the Cauchy-Schwarz 
inequality
$|Q(\gamma,c)|^2 \leq sQ(c,c)$. If $c$ is the Chern class of an almost 
complex structure
$J$, this implies $|Q(\gamma,c)|^2 \leq s(2\chi+3b_2)$. 
If $\Sigma$ is a representative of
$\gamma$ that is pseudoholomorphic with respect to $J$, we therefore have
$(2-2g+s)^2 \leq s(2\chi+3b_2)$. Note that the number on the right side 
of this inequality
must be non-negative, otherwise there is no almost complex structure on $X$ 
at all. 
Hence we see that there is an upper bound for the genus of pseudoholomorphic 
curves
representing $\gamma$ that does not depend on $J$, therefore a 
representative with large
genus provides the required example. A similar argument works if $Q$ is 
negative definite.
\end{proof}

\begin{proof} of Proposition \ref{propcountex}:
If $Q$ is definite, then the assertion of the
Proposition is covered by Proposition \ref{prop2}. If both $b^+$ and $b^-$
equal 1, then $Q$ must be the intersection form of $S^ 2 \times S^2$ or
of $\CP^2 \# \cpq$, hence the
intersection form of a rational surface. The same is true if $b^+=1$ and 
$Q$ is odd, all these cases are covered by Proposition \ref{prop1}.  
So the last case that is not covered by any of the preceeding examples is
the case $b^-=1$ and $b^ + \geq 2$. 
In this case, choose a 4-manifold $X'$ with intersection form $Q$. 
Let $h(X')=2 \chi(X') + 
3 \tau(X')$. 
Choose a class $\gamma$ in $H^2(X';\Z)$ with divisibility 1 and 
$-s^2:=Q(\gamma,\gamma) < 0$. 
Consider the lattice defined
by the integral cohomology in the semi-euclidean vector space $H^2(X;\R)$. 
Choose a basis $e_1, 
\dots ,e_n$ of this
vector space such that with respect to this basis, $Q$ is given by the matrix
$$
Q=\begin{pmatrix} -1 \\ & 1 \\ & & 1 \\ & & & \cdot \\ & & & & \cdot 
\\ & & & & & 1 
\end{pmatrix}
$$
where all other entries are zero, and such that $\gamma = s e_1$ 
(more precisely, 
the free part of $\gamma$).
Now choose a surface $\Sigma$ representing $\gamma$ with arbitrary genus $g$. 
Let $X=X' \# k (S^ 1 \times S^ 3)$,
where $k$ is chosen large, such that $h(X) = 2\chi(X) + 3 \tau(X) 
< -(\frac{2-2g}{-s}+s)^2$ 
and $b_1(X) + b_2^+(X) 
\equiv 1 \mod 2$ (note that attaching a copy
of $S^ 1 \times S^ 3$ decreases the Euler characteristic $\chi(X)$ by 2 
without changing the second 
cohomology and the intersection form).
If now $J$ would be any almost complex structure on $X$ such that 
$\Sigma$ is pseudoholomorphic 
with respect to $J$,
and $c=\sum_i c_i e_i \in H^ 2(X;\R)$ its (real) first Chern class  , 
then the adjunction 
equality would imply
$c \cdot \Sigma = -sc_1 = 2-2g-s^2$, 
hence $$
Q(c,c) = -c_1^ 2 + \sum _{i \neq 1} c_i^ 2 \geq -c_1^ 2 = 
-(\frac{2-2g}{-s}+s)^ 2 > h(X),
$$
in contradiction to $c^2 = 2 \chi(X) + 3 \tau(X) = h(X)$. This proves 
that $\Sigma$ can not be pseudoholomorphic for
any $J$, although, since the divisibility of $\gamma$ is one, the 
condition $g \equiv k(\gamma) \mod d$ is fulfilled.
\end{proof}

\bigskip

Finally I would like to thank my advisor, Dieter Kotschick, for his
support and the helpfull discussions during the work on this paper.

\end{document}